# MIXED-RATES ASYMPTOTICS

By Peter Radchenko

*University of Southern California*

A general method is presented for deriving the limiting behavior of estimators that are defined as the values of parameters optimizing an empirical criterion function. The asymptotic behavior of such estimators is typically deduced from uniform limit theorems for rescaled and reparametrized criterion functions. The new method can handle cases where the standard approach does not yield the complete limiting behavior of the estimator. The asymptotic analysis depends on a decomposition of criterion functions into sums of components with different rescalings. The method is explained by examples from Lasso-type estimation, $k$-means clustering, Shorth estimation and partial linear models.

**1. Introduction.** Consider an estimator $(a_n, b_n)$ that in some sense optimizes a random criterion function $G_n(a, b)$ over an open subset of $\mathbb{R}^{d_1} \times \mathbb{R}^{d_2}$. Two types of *mixed-rates* asymptotic behavior can occur and often occur simultaneously. First, the components $a_n$ and $b_n$ of the estimator may converge at different rates. Second, the criterion function itself may have important components settling down at different rates. The new method presented in this paper can handle both types of mixed-rates behavior.

Deriving the asymptotics of an estimator can be viewed as a three step procedure: proving consistency, establishing the rate of convergence and deriving the limiting distribution. This paper concentrates only on the last two steps. The limiting distribution is typically derived via a uniform limit theorem for the rescaled and reparametrized criterion functions. Suppose that the rates of convergence for the two components of the estimator have been established: $q_n^{-1}\|a_n - a_0\| \vee r_n^{-1}\|b_n - b_0\| = O_p(1)$ for some fixed parameter value $(a_0, b_0)$. Consider *localized criterion functions* of the form

$$H_n(s,t) := G_n(a_0 + q_n s, b_0 + r_n t) - G_n(a_0, b_0).$$

If, after appropriate rescaling, random functions $H_n(s,t)$ settle down to a "nice" stochastic process, the convergence in distribution of vectors $(s_n, t_n) := (q_n^{-1}[a_n - a_0], r_n^{-1}[b_n - b_0])$ to the corresponding optimizer of the limit process may follow from a continuous mapping type of argument. Theorem 3.2.2 of van der Vaart and Wellner [21] makes this argument precise for estimators defined by maximization. The above approach is standard when the rates $r_n$ and $q_n$ are the same, and it can work in some mixed-rates cases, such as the change-point problem (see, e.g., the section on nonregular examples in Kosorok [8]). Other mixed-rates examples where this argument succeeds can be found in Rotnitzky, Cox, Bottai and Robins [16], Pollard and Radchenko [12] and Andrews [1].

Many mixed-rates problems cannot be completely handled by the above approach. In the examples considered in this paper, the localized criterion function has the form

$$H_n(s,t) = \alpha_n f_n(s) + \beta_n g_n(s,t),$$

where $\beta_n = o(\alpha_n)$, the random function $f_n(s)$ settles down to a stochastic process $f(s)$, and $g_n$ is stochastically bounded. Because the limit of $\alpha_n^{-1} H_n(s,t)$ is a stochastic process indexed only by $s$, the standard approach fails to establish the limiting distribution of the component $t_n$. However, if random function $g_n(s,t)$ settles down to a stochastic process $g(s,t)$, a two-step continuous mapping argument can be used to establish the distributional limit of the vector $(s_n, t_n)$. This general idea is made rigorous by Theorem 1 in Section 2.

Another challenging problem is deriving the *correct* rates of convergence for the two components of the estimator. Standard methods represent the centered criterion function $G_n(a,b) - G_n(a_0, b_0)$ as a sum of a positive deterministic function and a random one, whose rates of growth around the value $(a_0, b_0)$ can be controlled (the deterministic function is typically approximated by a quadratic, and the random function is often approximately linear). Balancing out the two terms produces the rate of convergence: see, for example, Theorem 3.2.5 and Theorem 3.2.16 in van der Vaart and Wellner [21]. When $a_n$ and $b_n$ converge at different rates, this approach yields the "correct" rate only for the slower converging component. A reparametrization of the problem can sometimes be applied beforehand to sidestep this issue (for interesting examples, see the references at the end of the paragraph on the standard method for deriving the limiting distribution). Unfortunately, such a trick is not available in general, and a more careful treatment of the criterion function is required. To derive the rate for the faster converging component, say, $b_n$, Theorem 2 in Section 3 balances out the terms in a similar, but typically a more complicated, representation for the function $b \mapsto [G_n(a_n, b) - G_n(a_n, b_0)]$.

Section 4 is devoted to mixed-rates problems that arise in M-estimation. Consider a collection of functions $g_\theta(x)$ and an empirical measure $P_n$, corresponding to independent observations coming from a distribution $P$. Define the estimator $\theta_n$ as the minimizer of the criterion function $G_n(\theta) = P_n g_\theta$, and suppose that function $G(\theta) = \int g_\theta \, dP$ is minimized by $\theta_0$. The stochastic bound $\|\theta_n - \theta_0\| = o_p(1)$ usually follows from a uniform law of large numbers, and the central limit theorem for the estimator is typically derived from a quadratic approximation of the form

$$G_n(\theta) - G_n(\theta_0) \approx (\theta - \theta_0)' G''(\theta_0)(\theta - \theta_0) + n^{-1/2}(\theta - \theta_0)' Z_n,$$

under the regularity assumption that matrix $G''(\theta_0)$ is a positive definite matrix. If this regularity assumption breaks down and $G''(\theta_0)$ is singular, the approximation has to be carried out to higher order terms, which typically leads to mixed-rates situations that standard methods cannot handle. Theorem 3 covers exactly such cases. The form of the approximation to function $G(\theta)$ near $\theta_0$ determines the rates of convergence and the main features of the limiting behavior of the components of the estimator. Various remainder terms are handled by simple conditions imposed on functions $g_\theta$.

Mixed-rates behavior naturally arises in the estimation of semiparametric models. Most of the results in this paper do not directly apply to such problems, but, as the example in Section 8 demonstrates, some of the methods and ideas can be carried over.

For the simplicity of the presentation, the estimators and the criterion functions considered in this paper have at most two components converging at different rates. All the results can be easily extended to cover cases of more than two mixed-rates components.

This paper is organized as follows. Sections 2, 3 and 4 contain the general mixed-rates asymptotics results, namely, the limiting distribution theorem, the rates of convergence theorem and the M-estimation theorem. Proofs of these theorems are confined to Section 9. Sections 5, 6 and 7 contain applications of the general results to particular problems in Lasso-type estimation, shorth estimation and $k$-means clustering. Section 8 discusses a semiparametric example.

The abbreviation $Qf = \int f \, dQ$ is used throughout the paper for a given measurable function $f$ and a signed measure $Q$. In particular, given independent observations $X_i$ coming from a distribution $P$, let $P_n f$ denote $\sum_{i \le n} f(X_i)/n$ and define the empirical process $\nu_n$ on a class of functions $f$ by

$$f \mapsto \nu_n f = n^{1/2}(P_n - P)f = n^{-1/2} \sum_{i=1}^{n} [f(X_i) - Pf].$$

Write $\|\cdot\|_2$ for the $L_2(P)$ norm and say that a function $f$ is square-integrable if $\|f\|_2 < \infty$. Interpret $f(\theta) \gtrsim g(\theta)$ to mean that there exists a positive constant $c_0$ such that $f(\theta) \ge c_0 g(\theta)$ for all $\theta$ in a sufficiently small neighborhood of the origin. Analogously, interpret $\alpha_n \gtrsim \beta_n$ to mean $\alpha_n \geq c_0 \beta_n$ for all sufficiently large $n$.

**2. Limiting distribution.** Let the estimator $(a_n, b_n)$ converge in probability to a fixed parameter value $(a_0, b_0)$. Suppose that the rates of convergence $q_n$ and $r_n$ have been established for the components $a_n$ and $b_n$, respectively. Vector $(s_n, t_n) := (q_n^{-1}[a_n - a_0], r_n^{-1}[b_n - b_0])$ optimizes the localized criterion function $H_n(s, t)$ and satisfies the tightness condition $\|(s_n, t_n)\| = O_p^*(1)$. Focus on deriving the limiting distribution of $(s_n, t_n)$ when it is defined by minimization.

To avoid some measurability issues by allowing nonmeasurable maps, convergence in distribution (denoted by "$\rightsquigarrow$") is understood in the sense of Hoffmann–Jørgensen. An exposition of this general concept can be found in the monographs of Dudley [3] and van der Vaart and Wellner [21]. Let $\mathcal{B}_{\mathrm{loc}}(\mathbb{R}^d)$ be the space of all locally bounded real functions on $\mathbb{R}^d$. Convergence of the random processes considered in the examples of this paper is handled by equipping $\mathcal{B}_{\mathrm{loc}}(\mathbb{R}^d)$ with the metric $\rho$ for the topology of uniform convergence on compacta:

$$\rho(g, h) = \sum_{k=1}^{\infty} 2^{-k} \min[1, \rho_k(g, h)] \qquad \text{where } \rho_k(g, h) = \sup_{\|t\| \leq k} |g(t) - h(t)|.$$

In Theorem 1, convergence of the components of the criterion function should be understood with respect to this metric. The following continuity property of the arg min functional with respect to $\rho$ simplifies the statement of the theorem. Let $x^*$ be the *clean* minimum of a function $h$ in the sense that the strict inequality $h(x^*) < \inf_{\varepsilon \leq |x - x^*| \leq r} h(x)$ is satisfied for all positive $r$ and $\varepsilon$. Then

(1) $\qquad \rho(h_n, h) \to 0 \quad \text{implies} \quad \underset{|x - x^*| \leq r}{\arg\min}\, h_n(x) \to x^* \qquad \text{for each } r > 0.$

Note that the unique minimum of a continuous function is also its clean minimum over each large enough ball. In fact, lower semicontinuity of the function is sufficient. The proof of Theorem 1 would remain valid if $\mathcal{B}_{\mathrm{loc}}(\mathbb{R}^d)$ were equipped with a different metric $d$, as long as assumption (1) were imposed explicitly and formulated in terms of $d$.

The following result is stated in the cleanest form that covers the examples considered in the paper, thus, some of its conditions can be relaxed. See Remark for the alternative to the continuity assumption placed on the sample paths of the limit process $(f, g)$. Also note that the sample path properties required of the limit process need to hold only almost surely.

THEOREM 1. *Let $H_n$ be random criterion functions on $\mathbb{R}^{d_1} \times \mathbb{R}^{d_2}$ and let $(s_n, t_n)$ be random vectors in $\mathbb{R}^{d_1} \times \mathbb{R}^{d_2}$. Suppose that the following conditions are satisfied:*

  (i) $H_n(s,t) = \alpha_n f_n(s) + \beta_n g_n(s,t)$, where $f_n$ and $g_n$ are random functions on $\mathbb{R}^{d_1}$ and $\mathbb{R}^{d_1} \times \mathbb{R}^{d_2}$ respectively, while $\alpha_n$ and $\beta_n$ are positive numbers with $\beta_n = o(\alpha_n)$;
  (ii) $(f_n, g_n) \rightsquigarrow (f, g)$ and the limit process has continuous sample paths;
  (iii) $H_n(s_n, t_n) \leq \inf_{s,t} H_n(s,t) + o_p^*(\beta_n)$ and
  (iv) $\|(s_n, t_n)\| = O_p^*(1)$.

Assume that the sample paths of $f(\cdot)$ possess a unique minimum at a (random) point $s^*$ and the sample paths of $g(s^*, \cdot)$ possess a unique minimum at $t^*$. Then $(s_n, t_n) \rightsquigarrow (s^*, t^*)$.

REMARK. The assumptions on the sample paths of the limit process $(f, g)$ can be relaxed as follows. Assume that $s^*$ and $t^*$ are measurable random points such that for almost all sample paths of the limit process:

  (a) $s^*$ is the "clean" minimum of $f(\cdot)$,
  (b) $t^*$ is the "clean" minimum of $g(s^*, \cdot)$ and
  (c) for each ball $B$, the set of functions $\{g(\cdot, t) : t \in B\}$ is equicontinuous.

Theorem 1 can be generalized to cover cases where the optimizer is not defined by minimization or maximization. Suppose that vectors $(s_n, t_n)$ satisfy equalities $s_n = \Psi[H_n(\cdot, t_n)]$ and $t_n = \Phi[H_n(s_n, \cdot)]$ for certain maps $\Psi$ and $\Phi$. Assume that these maps are invariant to multiplications by positive constants and that $\Phi$ is also invariant to translations. If, in addition, each map satisfies assumption (1) with the proper replacement for the arg min, the proof of Theorem 1 still goes through. For a rigorous account of this fact, see Theorem 1 in Radchenko [15].

**3. Rates of convergence.** Consider two-component estimators $(a_n, b_n)$ that are defined by minimizing random criterion functions $G_n(a, b)$. The following lemma uses an approximation to the criterion function to establish the rate of convergence of the slower converging component $a_n$ and makes an initial guess at the rate of convergence of the component $b_n$. This guess is not quite correct, but it provides an improvement over existing results, which establish one convergence rate for the whole long vector $(a_n, b_n)$. Lemma 1 requires a particular representation for the criterion function. In many standard asymptotic problems, this representation is satisfied with the term $M_n(a, b)$ bounded below by a nonsingular quadratic, and the term $N_n(a, b)$ of the order $O_p(n^{-1/2}\|(a, b)\|)$, which yields the usual $n^{-1/2}$ rate of convergence. The lemma handles cases that are more general.

LEMMA 1. *Suppose that inequalities $G_n(a_n, b_n) \leq G_n(0, 0)$ hold together with the stochastic bound $\|(a_n, b_n)\| = o_p^*(1)$. Let $\alpha$ and $\beta$ be positive numbers*

satisfying $\alpha \geq \beta$, and let $\{\gamma_1, \ldots, \gamma_p, \eta_1, \ldots, \eta_p\}$ be a collection of nonnegative numbers satisfying $\gamma_i < \alpha$ for all $i \in \{1, \ldots, p\}$. Suppose that criterion functions $G_n$ satisfy a representation

$$G_n(a,b) - G_n(0,0) = M_n(a,b) - N_n(a,b),$$

such that

$$M_n(a_n, b_n) \gtrsim \|a_n\|^\alpha + \|b_n\|^\beta \qquad \text{with inner probability tending to one, and}$$

$$[N_n(a_n, b_n)]^+ = O_p^*\left(\sum_{i \leq p} n^{-\eta_i} \|(a_n, b_n)\|^{\gamma_i}\right).$$

Define $\tau_a = \min_{i \leq p}(\frac{\eta_i}{\alpha - \gamma_i})$. Then $\|a_n\| = O_p^*(n^{-\tau_a})$ and $\|b_n\| = O_p^*(n^{-\alpha \tau_a/\beta})$.

Once the convergence rate of $a_n$ is established, it becomes reasonable to fix $a = a_n$ and consider the function $b \mapsto G_n(a_n, b)$. Existing results do not necessarily yield the convergence rate of the minimizer of this function. The point of difficulty is that the leading terms in the approximation to this function near its minimum are more complex than the ones that appear in the standard asymptotics. The following theorem can handle such cases but it requires a more refined approximation to the criterion function. One may want to use the help of Lemma 1 to obtain such an approximation (see, e.g., the proof of Theorem 3), and then apply Theorem 2 to derive the "correct" convergence rate of $b_n$. Note that Theorem 2 places no assumptions at all on the space containing the $a$-component.

THEOREM 2. *Let $G_n(a,b)$ be a function of two components, where the first component belongs to an abstract set, and the second belongs to a Euclidean space. Suppose that inequalities $G_n(a_n, b_n) \leq G_n(a_n, 0)$ hold together with the stochastic bound $\|b_n\| = o_p^*(1)$. Let $\beta$ be positive and let $\{\alpha_1, \ldots, \alpha_p, \beta_1, \ldots, \beta_p\}$ be a collection of nonnegative numbers satisfying $\beta_i < \beta$ for all $i \in \{1, \ldots, p\}$. Assume that $G_n$ satisfies a representation*

(2) $$G_n(a,b) - G_n(a,0) = M_n(a,b) - N_n(a,b),$$

*such that*

$$M_n(a_n, b_n) \gtrsim \|b_n\|^\beta \qquad \text{with inner probability tending to one, and}$$

$$[N_n(a_n, b_n)]^+ = O_p^*\left(\sum_{i \leq p} n^{-\alpha_i} \|b_n\|^{\beta_i}\right).$$

*Then $\|b_n\| = O_p^*(n^{-\tau_b})$ for $\tau_b = \min_{i \leq p}\{\frac{\alpha_i}{\beta - \beta_i}\}$. If $[N_n]^+ \equiv 0$ then $\mathbb{P}_*\{b_n = 0\} \to 1$.*

**4. M-estimators.** The following definition introduces notation that is used in the statement of Theorem 3. This notation simplifies the work with polynomials that are homogeneous functions of the elements of vector $(a,b)$ and the absolute values of the elements of vector $(a,b)$.

DEFINITION 1. Let $\psi$ be a real valued function on $\mathbb{R}^d$ and let $\gamma$ be a positive constant. Say that $\psi \in H_1^+(\gamma)$ if $\psi(\lambda\theta) = \lambda^\gamma \psi(\theta)$ for all $\lambda \geq 0$ and $\psi(\theta) > 0$ for all $\theta \neq 0$.

Let $\phi$ be a real valued function on $\mathbb{R}^{d_1} \times \mathbb{R}^{d_2}$ and let $\alpha$ and $\beta$ be some positive constants. Say that $\phi \in H_2^{(-)}(\alpha, \beta)$ if $\phi(\lambda_1 a, \lambda_2 b) = (\lambda_1)^\alpha (\lambda_2)^\beta \phi(a,b)$ for all nonnegative $\lambda_1$ and $\lambda_2$, while function $\phi$ assumes at least some negative values.

REMARK. For each continuous function $\psi(\theta)$ in the class $H_1^+(\gamma)$, there exist positive constants $c_1$ and $c_2$ such that $c_1 \|\theta\|^\gamma \leq \psi(\theta) \leq c_2 \|\theta\|^\gamma$.

Suppose that $X_1, X_2, \ldots, X_n$ are independent observations in $\mathbb{R}^k$ coming from a distribution $P$ and write $P_n$ for the corresponding empirical distribution. Suppose that $A$ is an open subset of $\mathbb{R}^{d_1} \times \mathbb{R}^{d_2}$ and let $\{g_{a,b}(x) : (a,b) \in A\}$ be a collection of real valued $P$-integrable functions on $\mathbb{R}^k$. Assume that this collection of functions is centered to satisfy $g_{0,0} \equiv 0$. Suppose that vectors $(a_n, b_n)$ minimize over $A$ the random criterion functions $G_n(a,b) = P_n g_{a,b}$ and let $(0,0)$ be the corresponding minimizer of the population analog $G(a,b) = P g_{a,b}$. The following theorem derives the asymptotics of $(a_n, b_n)$ in the challenging case of the singular second derivative matrix $G''(0,0)$.

THEOREM 3. *Let $\{\alpha, \beta, \gamma_1, \ldots, \gamma_p, \eta_1, \ldots, \eta_p\}$ be a collection of positive numbers. Assume that $\alpha > \beta > 1$ and $\beta > \eta_j$ for $1 \leq j \leq p$. Suppose that there exist continuous functions $\psi_1(a) \in H_1^+(\alpha), \psi_2(b) \in H_1^+(\beta)$ and $\phi_i(a,b) \in H_2^{(-)}(\gamma_i, \eta_i)$ for $1 \leq i \leq p$, such that near the origin the population criterion function satisfies the following conditions:*

(i) $G(a,b) \gtrsim \|a\|^\alpha + \|b\|^\beta$,
(ii) $G(a,0) = \psi_1(a) + o(\|a\|^\alpha)$ and
(iii) $G(a,b) = G(a,0) + \psi_2(b)[1 + o(1)] + \sum_{i=1}^p \phi_i(a,b)[1 + o(1)] + o(\sum_{i=1}^\alpha \|a\|^{\alpha-i}\|b\|^i)$.

*Let $\tau_a = \frac{1}{2(\alpha-1)}$, $\lambda_0 = \frac{1}{2(\beta-1)}$, $\lambda_i = \frac{\tau_a \gamma_i}{\beta - \eta_i}$ for $1 \leq i \leq p$, and define $\tau_b = \min_{0 \leq i \leq p}[\lambda_i]$. Suppose there exist on $\mathbb{R}^k$ five square integrable functions, $\Delta_1$ (taking values in $\mathbb{R}^{d_1}$), $\Delta_2$ (taking values in $\mathbb{R}^{d_2}$) and real valued $r_{a,b}, s_{a,b}$ and $l_{a,b}$, such that:*

(iv) $g_{a,b}(x) = a'\Delta_1(x) + b'\Delta_2(x) + \|(a,b)\| r_{a,b}(x)$;
(v) $g_{a,b}(x) - g_{a,0}(x) - b'\Delta_2(x) = l_{a,b}(x) + \|b\| s_{a,b}(x)$;

(vi) $\sup_{\|(a,b)\|\leq \delta_n} |\nu_n r_{a,b}| = o_p(1)$ and $\sup_{\|(a,b)\|\leq \delta_n} |\nu_n s_{a,b}| = o_p(1)$ for all $\delta_n \to 0$;

(vii) $\sup_{\|a\|\leq \delta_n, \|b\|\leq \varepsilon_n} |\nu_n l_{a,b}| = o_p(n^{-\beta\tau_b+1/2})$ for all $\delta_n = O(n^{-\tau_a}), \varepsilon_n = O(n^{-\alpha\tau_a/\beta})$.

Assume that $\|(a_n, b_n)\| = o_p(1)$. If $\alpha\tau_a = \beta\tau_b$, then

$$(n^{\tau_a}a_n, n^{\tau_b}b_n) \rightsquigarrow \arg\min_{s,t}\left[\psi_1(s) + s'Z_1 + \psi_2(t) + 1\{\lambda_0 = \tau_b\}t'Z_2 + \sum_{i=1}^{p} 1\{\lambda_i = \tau_b\}\phi_i(s,t)\right];$$

otherwise $(n^{\tau_a}a_n, n^{\tau_b}b_n) \rightsquigarrow (s^*, t^*)$, where

$$s^* = \arg\min_{s}[\psi_1(s) + s'Z_1],$$

$$t^* = \arg\min_{t}\left[\psi_2(t) + 1\{\lambda_0 = \tau_b\}t'Z_2 + \sum_{i=1}^{p} 1\{\lambda_i = \tau_b\}\phi_i(s^*,t)\right].$$

Here $(Z_1, Z_2)$ is a mean zero Gaussian vector with covariance matrix $P(\Delta_1, \Delta_2) \times (\Delta_1, \Delta_2)'$.

Note that a stochastic process $\nu_n f_{a,b}$ necessarily satisfies the uniform stochastic bound required in condition (3) of the above theorem (cf. *asymptotic equicontinuity* defined in van der Vaart [20]) if functions $f_{a,b}$ form a Donsker class and $\|f_{a,b}\|_2 \to 0$ as $\|(a,b)\| \to 0$. Simple ways of checking that a class of functions is Donsker are given, for example, in van der Vaart's Theorem 19.5 and Theorem 19.14.

To illustrate the variety of asymptotic results produced by Theorem 3, consider some simple approximations to the function $G$, which has a singular second derivative at the origin, where its minimum is located. Let $(a,b) \in \mathbb{R}^2$ and consider the case $G(a,b) \approx a^4 + b^2$. Theorem 3 yields $(n^{1/6}a_n, n^{1/2}b_n) \rightsquigarrow (\arg\min_s[s^4 + sZ_1], \arg\min_t[t^2 + tZ_2])$ if the conditions (iv)–(vii) are satisfied. Here $(Z_1, Z_2)$ is a mean zero Gaussian vector. Now consider the case $G(a,b) \approx a^4 + b^2 + a^2b$. Under the same assumptions, the theorem yields $(n^{1/6}a_n, n^{1/3}b_n) \rightsquigarrow (\arg\min_{s,t}[s^4 + sZ_1 + t^2 + s^2t])$. If the approximation is $G(a,b) \approx a^4 + b^2 + a^3b$, the corresponding result is $(n^{1/6}a_n, n^{1/2}b_n) \rightsquigarrow (s^*, t^*)$ with $s^* = \arg\min_s[s^4 + sZ_1]$ and $t^* = \arg\min_t[t^2 + (s^*)^3t + tZ_2])$. Note that Theorem 3 does not attempt to cover every conceivable approximation to $G(a,b)$, as the statement of the result would become too long and complicated, but each such situation can be handled with only minor modifications to the proof of the theorem.

**5. Example: Lasso-type estimators.** Assume that the observed variables $Y_i$ satisfy the linear model

$$Y_i = x_i'\beta + \varepsilon_i, \qquad i = 1, \ldots, n.$$

The errors $\varepsilon_i$ are independent and identically distributed random variables that have mean zero and variance $\sigma^2$. The parameter $\beta$ is a vector in $\mathbb{R}^d$ that needs to be estimated. The covariates $x_i$ are fixed and centered, and the matrix $C_n = \frac{1}{n}\sum_{i=1}^n x_i x_i'$ is nonsingular.

Suppose $\lambda_n$ and $\gamma$ are positive real numbers. Define the "Lasso-type" estimator $\beta_n$ as the minimizer of the penalized least-squares criterion,

$$W_n(\alpha) = \sum_{i=1}^n (Y_i - x_i'\alpha)^2 + \lambda_n \sum_{j=1}^d |\alpha_j|^\gamma,$$

over all vectors $\alpha = (\alpha_1, \ldots, \alpha_d)'$. In the particular cases of $\gamma = 1$ and $\gamma = 2$, this estimator corresponds, respectively, to the "Lasso" of Tibshirani [18] and the ridge regression. For general $\gamma$, such estimators were introduced by Frank and Friedman [4]. The limiting behavior of the estimator $\beta_n$ was described by Knight and Fu [7] under certain conditions on the growth rate of the weighting sequence $\{\lambda_n\}$.

Assume that the design satisfies the following regularity conditions:

(i) matrixes $C_n$ converge to a fixed matrix $C$;
(ii) as $n$ tends to infinity, $n^{-1}\max_{i \le n}(x_i'x_i)$ converges to zero.

In the case of the nonsingular matrix $C$, Knight and Fu derived the $\sqrt{n}$-asymptotics for $\beta_n$ after setting the growth rate for the weighting sequence $\{\lambda_n\}$. They required that, for some nonnegative constant $\lambda_0$,

(3) $$\lambda_n / n^{\min(1/2, \gamma/2)} \to \lambda_0.$$

Note that when $\lambda_0 = 0$, the penalty contribution is asymptotically negligible and the limiting behavior of the estimator $\beta_n$ is the same as that of the usual least-squares estimator.

To derive the asymptotics of $\beta_n$, Knight and Fu used a standard approach that is based on rescaling the parameters at the same rate and applying a continuous mapping type of argument. When vector $\beta$ has a zero component, $\gamma < 1$, and $\lambda_n$ grows faster than the rate given in (3), this approach fails to deliver the complete asymptotics. For concreteness, consider the case $d = 2, \beta = (1, 0)', \gamma = 1/2$, and set $\lambda_n = \lambda_0 n^{1/2}$ for some positive constant $\lambda_0$. The standard approach establishes the asymptotics of the first component of $\beta_n$, but only yields the $o_p(n^{-1/2})$ stochastic order for the second component of the estimator (see Knight and Fu [7], page 1361). The techniques developed in Section 3 are applied below to show that the second component is in fact exactly zero with probability tending to one.

Because $C$ is nonsingular and $\lambda_n = o(n)$, the estimator $\beta_n$ is consistent (see Theorem 1 of Knight and Fu [7]). The proof is based on the fact that, for each fixed $\alpha$, the penalty part of the criterion function $W_n(\alpha)$ is asymptotically negligible compared to the least-squares part. Focus on vectors $\alpha$ that are near the true parameter $\beta$, and write $\alpha = \beta + (a,b)'$. Express the penalized criterion function in terms of $a$ and $b$. Denote $n^{-1}[W_n(\alpha) - W_n(\beta)]$ by $G_n(a,b)$, and let $Z_n$ stand for $n^{-1/2}\sum_{i=1}^n \varepsilon_i x_i$. The regularity conditions on the design guarantee that the sequence of random vectors $Z_n$ has a limiting Gaussian distribution with mean zero and covariance $\sigma^2 C$. As $a$ and $b$ tend to zero,

$$G_n(a,b) = (a,b)C_n(a,b)' - 2n^{-1/2}(a,b)Z_n + \frac{\lambda_0}{2}n^{-1/2}a[1+o(1)] + \lambda_0 n^{-1/2}|b|^{1/2}.$$

The $o(1)$ terms come from the Taylor expansion of $|1+a|^{1/2}$ near $a=0$. Function $G_n$ is minimized by the vector $(a_n, b_n)'$ that is defined as the difference between $\beta_n$ and $\beta$.

Define $M_n(a,b)$ to be $(a,b)C_n(a,b)' + \lambda_0 n^{-1/2}|b|^{1/2}$ and let $N_n$ equal $G_n - M_n$. Note that for all $n$ large enough, the eigen values of the sequence of matrixes $C_n$ are bounded away from zero. Apply Lemma 1 from Section 3 and conclude that $\|(a,b)\| = O_p(n^{-1/2})$.

Let $v_n$ denote the bottom right element of the matrix $C_n$. Observe that

$$G_n(a_n, b_n) - G_n(a_n, 0) = v_n b_n^2 + O_p(n^{-1/2}|b_n|) + \lambda_0 n^{-1/2}|b_n|^{1/2}.$$

Note that $\lambda_0$ and $v_n$ are positive and $v_n$ is bounded away from zero for all sufficiently large $n$. Deduce that, with probability tending to one, the right-hand side of the above display is bounded below by $cb_n^2$ for some positive $c$. Apply Theorem 2 with $N_n \equiv 0$ and conclude that $\mathbb{P}\{b_n = 0\} \to 1$.

More examples of mixed-rates behavior in Lasso-type estimation can be found in Radchenko [14].

**6. Example: Shorth.** Assume that the observations are independently sampled from a distribution $P$ on the real line and let $[m_n - r_n, m_n + r_n]$ be the shortest interval that contains at least half of the first $n$ observations. The shorth estimator is defined as the average over such an interval, but the goal of this section is the limiting behavior of $m_n$ and $r_n$. Grübel [5] derived the root-$n$ asymptotics for $r_n$ and Kim and Pollard [6] derived the cube root asymptotics for $m_n$. The methods of the present paper allow one to establish the joint limiting behavior of $(m_n, r_n)$ using a simple approximation to the criterion function.

Denote by $\mu$ and $\rho$ the population solution, in other words, let $[\mu - \rho, \mu + \rho]$ be the shortest interval to which $P$ assigns at least half the probability.

Assume that the population solution is unique and let $P$ have a bounded density $f$ that is differentiable at the endpoints $\mu \pm \rho$. Define the criterion function $V_n$ by $V_n(\varepsilon, \delta) = P_n[(\mu + \varepsilon) - (\rho + \delta), (\mu + \varepsilon) + (\rho + \delta)] - 1/2$, and let $V(\varepsilon, \delta)$ denote the population analog obtained by replacing $P_n$ with $P$. Observe that $V(0,0) = 0$ and write out the Taylor expansion for function $V$ near the origin:

$$(4) \qquad V(\varepsilon, \delta) = c_1 \delta + c_2 \varepsilon^2 + c_3 \varepsilon \delta + c_4 \delta^2 + o(\varepsilon^2 + \delta^2),$$

where the coefficients are $c_1 = f(\mu - \rho) + f(\mu + \rho), c_2 = c_4 = [f'(\mu + \rho) - f'(\mu - \rho)]/2$ and $c_3 = f'(\mu + \rho) + f'(\mu - \rho)$. The coefficient of the linear term in $\varepsilon$ equals zero because the function $V(\varepsilon, 0)$ is maximized at $\varepsilon = 0$. This forces the equality $f(\mu + \rho) = f(\mu - \rho)$. By the same reasoning, coefficient $c_2$ must be nonpositive. Assume $c_2 < 0$ and $c_1 > 0$ for regularity.

Recall the bound $\sup_{m,r} |P_n[m - r, m + r] - P[m - r, m + r]| = O_p(n^{-1/2})$ from the standard empirical process theory. Denote this supremum by $\Delta_n$. Uniqueness of the population solution and regularity assumptions on the coefficients of the Taylor expansion (4) guarantee that there exists a positive constant $c$ such that, for all small enough positive $\delta$, inequality $\sup_m P[m - (\rho - \delta), m + (\rho - \delta)] < 1/2 - c\delta$ holds. Consequently,

$$\sup_m P_n[m - (\rho - \Delta_n/c), m + (\rho - \Delta_n/c)] < \Delta_n + 1/2 - c\Delta_n/c = 1/2,$$

and hence, $\delta_n \geq -\Delta_n/c$. Expansion (4) also implies existence of a positive constant $b$ such that $P[\mu - (\rho + \delta), \mu + (\rho + \delta)] \geq 1/2 + b\delta$ for all small enough positive $\delta$. Take $\delta = \Delta_n/b$ and deduce that $P_n[\mu - (\rho + \Delta_n/b), \mu + (\rho + \Delta_n/b)] \geq 1/2$. Conclude that $\delta_n \leq \Delta_n/b$ and, hence, $\delta_n = O_p(n^{-1/2})$.

Note that function $V_n(\cdot, \delta_n)$ is maximized by $\varepsilon_n$ and function $V(\cdot, 0)$ has a clean maximum at zero. Uniform convergence in probability of $V_n(\cdot, \delta_n)$ to $V(\cdot, 0)$ implies $\varepsilon_n = o_p(1)$. Introduce functions $M_n(\varepsilon, \delta) = |c_2|\varepsilon^2/4$ and define functions $N_n$ by equalities

$$-[V_n(\varepsilon, \delta) - V_n(0, \delta)] = M_n(\varepsilon, \delta) - N_n(\varepsilon, \delta).$$

Note that when $\delta = \delta_n$, the expression on the left-hand side is minimized by $\varepsilon = \varepsilon_n$. Denote the difference between the indicator functions of intervals $[(\mu + \varepsilon) - (\rho + \delta), (\mu + \varepsilon) + (\rho + \delta)]$ and $[\mu - (\rho + \delta), \mu + (\rho + \delta)]$ by $J(\varepsilon, \delta)$. Observe that

$$V_n(\varepsilon, \delta) - V_n(0, \delta) = V(\varepsilon, \delta) - V(0, \delta) + (P_n - P)J(\varepsilon, \delta).$$

Recall that $c_2 < 0$ by the regularity assumptions placed on the coefficients of expansion (4), and use the Taylor expansion (4) to deduce a stochastic bound

$$(5) \quad N_n(\varepsilon_n, \delta_n) \leq (P_n - P)J(\varepsilon_n, \delta_n) - |c_2|\varepsilon_n^2/2 + O_p(n^{-1/2}|\varepsilon_n|) + O_p(n^{-1}).$$

Note that the collection of functions $J(\varepsilon, \delta)$ is a Vapnik–Cervonenkis class. For $R$ near zero, the envelope function $G_R = \sup_{\{\varepsilon^2 + \delta^2 \leq R^2\}} |J(\varepsilon, \delta)|$ is the indicator of the two intervals of total length bounded above by $4R$; boundedness of the density implies $PG_R^2 = O(R)$. Hence, the conditions of Lemma 4.1 of Kim and Pollard [6] are satisfied and, consequently, the bound $|(P_n - P)J(\varepsilon_n, \delta_n)| - c\varepsilon_n^2 \leq O_p(n^{-2/3})$ is valid for each positive $c$. It follows that

$$[(P_n - P)J(\varepsilon_n, \delta_n) - |c_2|\varepsilon_n^2/2]^+ = O_p(n^{-2/3}).$$

Combine this stochastic bound with bound (5) and deduce that

$$[N_n(\varepsilon_n, \delta_n)]^+ = O_p(n^{-1/2}|\varepsilon_n|) + O_p(n^{-2/3}).$$

An application of Theorem 2 yields $\varepsilon_n = O_p(n^{-1/3})$.

Set $I_{s,t} = 1[(\mu + n^{-1/3}t) - (\rho + n^{-1/2}s), (\mu + n^{-1/3}t) + (\rho + n^{-1/2}s)] - 1[\mu - \rho, \mu + \rho]$ and define the localized criterion functions $H_n(s, t) = V_n(n^{-1/3}t, n^{-1/2}s)$. Use the empirical process notation to write an approximation to $H_n(s, t)$ that holds uniformly on compacta:

$$
\begin{aligned}
H_n(s, t) = {} & n^{-1/2}\{c_1 s + \nu_n[\mu - \rho, \mu + \rho]\} \\
& + n^{-2/3}\{c_2 t + n^{1/6}\nu_n I_{s,t} + o(1)\}.
\end{aligned}
\tag{6}
$$

On each compact set, the stochastic processes $X_n(s, t) = n^{1/6}\nu_n I_{n^{1/6}s,t}$ converge in distribution to a tight Gaussian process by Theorem 19.28 of van der Vaart [20]. The conditions of the theorem are checked in van der Vaart's Example 19.29 for essentially the same process as $X_n$. Consequently, $X_n$ satisfies the asymptotic equicontinuity condition of van der Vaart's Theorem 18.14, and approximation $n^{1/6}\nu_n[I_{s,t} - \nu_n I_{0,t}] = o_p(1)$ holds uniformly over $s$ and $t$ in each given compact set. Note that

$$\lim_{n \to \infty} n^{1/3} P I_{0,t} I_{0,t'} = c_1 \min(|t|, |t'|)$$

and

$$\lim_{n \to \infty} n^{1/6} P I_{0,t} 1[\mu - \rho, \mu + \rho] = 0.$$

Write $f_n(s)$ and $g_n(s, t)$ for the two expressions in curly brackets that appear in representation (6). Let $B(t)$ be a two-sided Brownian motion and let $Z$ be an independent $N(0, 1/2)$ random variable. Conclude that

$$(f_n(s), g_n(s, t)) \rightsquigarrow (c_1 s - Z, c_2 t + \sqrt{c_1} B(t)).$$

Recall that the rescaled solution $(s_n, t_n) = (n^{1/2}\delta_n, n^{1/3}\varepsilon_n)$ is stochastically bounded and note the relationship

$$s_n = \inf\{s : H_n(s, t_n) \geq 0\} \quad \text{and} \quad t_n = \arg\max[H_n(s_n, \cdot)].$$

Consider the functional $\Psi:h \mapsto \inf\{s:h(s) \geq 0\}$ and note that the value $\Psi[H_n(\cdot,t)]$ is well defined and finite for each $t$. Also note that $\Psi$ is invariant to multiplications by positive constants. Apply Theorem 1 in Radchenko [15] and express the result in terms of the original variables:

$$(n^{1/2}[r_n - r], n^{1/3}[m_n - \mu]) \rightsquigarrow \left(Z/c_1, \arg\max_t [c_2 t + \sqrt{c_1} B(t)]\right).$$

Standard techniques fail to extract the limiting behavior of the estimator directly from approximation (6) because the first component of the approximation dominates the essential second component as $n$ tends to infinity.

**7. Example: $k$-means.** The $k$-means procedure divides observations $X_1, \ldots, X_n$ in $\mathbb{R}^d$ into $k$ sets by locating the cluster centers and then assigning each observation to the closest center. The set of cluster centers $C_n = \{c_{1n}, \ldots, c_{kn}\}$ is chosen to minimize

(7)  $$W_n(C) = n^{-1} \sum_{i \leq n} \min_{1 \leq j \leq k} \|x_i - c_j\|^2$$

as a function of sets $C = \{c_1, \ldots, c_k\}$ of $k$ not necessarily distinct points in $\mathbb{R}^d$. Assume that the observations are independent and come from a distribution $P$ on $\mathbb{R}^d$. Define the population criterion function, $W(C) = P \min_{1 \leq j \leq k} \|X_1 - c_j\|^2$, and let $C_0$ be a set that minimizes $W$. Note that if $P$ has a finite second moment and is not concentrated on fewer than $k$ points, then each set of optimal population centers has to contain exactly $k$ points. Under these conditions, and given that the set $C_0$ of optimal population centers is uniquely defined, Pollard [9] showed that the sets $C_n$ of optimal empirical centers are strongly consistent with respect to the Hausdorff metric.

In the example that follows, condition (vi) of Theorem 3 needs to be verified for classes of functions that possess the following simple property.

PROPERTY 1. *The class of functions $f_\theta(\cdot)$ satisfies the following conditions:*

(i) *the envelope function $F(\cdot)$ is square integrable with respect to $P$;*
(ii) *there exist positive integers $N$ and $d$ such that each $f_\theta$ can be represented as a sum of at most $N$ functions of the form $LQ$, where $L$ is a linear function and $Q$ is the indicator function of the intersection of at most $N$ half-spaces in $\mathbb{R}^d$;*
(iii) $\|f_\theta\|_2 \to 0$ *as $\theta \to 0$.*

The first two conditions imply that the class of functions $f_\theta$ is Donsker. This fact is proved on page 921 of Pollard [10], but it can also be easily deduced from the standard results on pages 274–276 of van der Vaart [20]. The third condition together with the Donsker property yield the required $\sup_{\|\theta\| \leq \delta_n} |\nu_n f_\theta| \to 0$ for each $\delta_n \to 0$.

The following is a two-dimensional extension of the example discussed in Section 7.1 of Radchenko [15]. Consider a distribution $Q$ on the plane $(x, y)$ that concentrates on the lines $\{x = 1\}$ and $\{x = -1\}$. Let $Q$ put probability one half on each line, and let the conditional distribution on each line be the double exponential. Write $Q$ as $P \times \mu$, where $P$ is the double exponential distribution and $\mu\{-1\} = \mu\{1\} = 1/2$.

There are two pairs of optimal population centers, $\{(-1,0),(1,0)\}$ and $\{(0,-1),(0,1)\}$; denote them by $C^v = \{c_1^v, c_2^v\}$ and $C^h = \{c_1^h, c_2^h\}$, respectively. The superscripts reflect either the vertical or the horizontal direction of the *split line*, which is defined as the common boundary for the two Voronoi half-planes generated by a given pair of centers. Let $C_n^v$ and $C_n^h$ minimize the criterion function (7) over two fixed nonoverlapping Hausdorff neighborhoods of the sets $C^v$ and $C^h$, respectively, and let $C_n$ be a global minimizer. A slight extension of Pollard's consistency result yields

$$C_n \in \{C_n^v, C_n^h\} \qquad \text{with probability tending to one,}$$

$$C_n^h \to C^h \quad \text{and} \quad C_n^v \to C^v \qquad \text{almost surely,}$$

where the set convergence is understood with respect to the Hausdorff metric. In fact, the probability with which $C_n$ chooses between the two configurations converges to a half. Near the set $C^h$, the population criterion function $W$ is approximated by a nonsingular quadratic. As a result, the solution $C_n^h$ settles down at the standard $n^{-1/2}$ rate and satisfies a central limit theorem. The remainder of the section is concerned with deriving the asymptotics of $C_n^v$, which is a challenging problem because the quadratic approximation to the population criterion function near the set $C^v$ is singular.

Suppose that $C = \{c_1, c_2\}$ is a candidate to minimize the criterion function (7) over a small Hausdorff neighborhood of the set $C^v$. Let $c_1 = (c_{1x}, c_{1y})$ be the point lying close to $c_2^v = (-1, 0)$ and let $c_2 = (c_{2x}, c_{2y})$ lie close to $c_2^v = (1, 0)$. Write $z$ to denote a point on the plane and let $(x, y)$ be the coordinate form of $z$. Introduce new variables by

$$\delta_s = \tfrac{1}{2}(c_{1x} + c_{2x}), \qquad \delta_d = 1 + \tfrac{1}{2}(c_{1x} - c_{2x}),$$
$$\varepsilon_s = \tfrac{1}{2}(c_{1y} - c_{2y}) \quad \text{and} \quad \varepsilon_d = \tfrac{1}{2}(c_{1y} + c_{2y}).$$

These variables contain the information on how far the centers in the set $C$ lie from the corresponding centers in $C^v$. Define $a = (\delta_s, \varepsilon_d)$ and $b = (\delta_d, \varepsilon_s)$. Let $g_{a,b}(z)$ be the squared distance from $z$ to the closest center in $C$, written in terms of $(a,b)$ and centered:

$$\begin{aligned}g_{a,b}(z) = &[(x+1-\delta_s-\delta_d)^2 + (y-\varepsilon_s-\varepsilon_d)^2] \\ &\wedge [(x-1-\delta_s+\delta_d)^2 + (y-\varepsilon_s+\varepsilon_d)^2] \\ &- \|z+(1,0)\|^2 \wedge \|z-(1,0)\|^2.\end{aligned}$$

Define criterion functions $G_n(a,b) = P_n g_{a,b}$ and $G(a,b) = P g_{a,b}$, and note that they are just the functions $W_n(C)$ and $W(C)$ centered at the set $C^v$ and written in terms of the new variables. Note that $G(a,b)$ is minimized at zero and the points $(a_n, b_n)$ that minimize $G_n(a,b)$ are of order $o_p(1)$ because of consistency.

The following approximation holds for $G$ near zero (see Section 2.3 of Radchenko [13]):

$$\begin{aligned}G(a,b) = &\tfrac{1}{6}(|\delta_s|+|\varepsilon_d|)^3 + \tfrac{1}{6}||\delta_s|-|\varepsilon_d||^3 \\ &+ \delta_d^2 + \varepsilon_s^2 + \delta_s^2\delta_d + 2\delta_s\varepsilon_d\varepsilon_s - \varepsilon_d^2\delta_d + O(\|(a,b)\|^4).\end{aligned}$$

Note that conditions (i), (ii) and (iii) of Theorem 3 are satisfied with $\alpha = 3, \beta = 2, p = 3$ and $(\gamma_i, \eta_i) = (2,1)$ for $1 \leq i \leq 3$. The corresponding homogeneous functions are

$$\psi_1(a) = \tfrac{1}{6}(|\delta_s|+|\varepsilon_d|)^3 + \tfrac{1}{6}||\delta_s|-|\varepsilon_d||^3, \qquad \psi_2(b) = \delta_d^2 + \varepsilon_s^2$$

and

$$\phi_1(a,b) = \delta_s^2\delta_d, \qquad \phi_2(a,b) = 2\delta_s\varepsilon_d\varepsilon_s, \qquad \phi_3(a,b) = -\varepsilon_d^2 d_d.$$

Take functions $\Delta_1(z)$ and $\Delta_2(z)$ in condition (iv) as the $L_2$ partial derivatives of $g_{a,b}(z)$ with respect to $a$ and $b$ at $(0,0)$. For example, let

$$b'\Delta_2(z) = -[2\delta_d(x+1) + 2\varepsilon_s y]H_-(z) + [2\delta_d(x-1) - 2\varepsilon_s y]H_+(z),$$

where $H_-(z)$ is the indicator function of the half-plane $\{z : x \leq 0\}$ and $H_+(z)$ is the indicator functions of the half-plane $\{z : x > 0\}$. Condition (vi) for the remainder functions $r_{a,b}(x)$ follows from a general result on $k$-means (see Pollard [10], Lemma B). The proof essentially consists of verifying that Property 1 holds for the class $\{r_{a,b}\}$.

The expression for $g_{a,b}(z)$ depends on the sign of $x$ and on which of the centers in the set $C$ lies closer to the point $z$. Let $D$ and $U$ be the $x$-coordinates of the crossing points of the split line corresponding to $C$ with the lines $\{y = -1\}$ and $\{y = 1\}$, respectively. Note that when $b = 0$, the values $D$ and $U$ are simply $\delta_s - \varepsilon_d$ and $\delta_s + \varepsilon_d$. Introduce functions

$$A(z) = 1\{|x| \leq |D|, xD > 0, y = -1\} + 1\{|x| \leq |U|, xU > 0, y = 1\}$$

and
$$A^0(z) = A(z) - 1\{|x| \leq |\delta_s - \varepsilon_d|, x(\delta_s - \varepsilon_d) > 0, y = -1\}$$
$$- 1\{|x| \leq |\delta_s + \varepsilon_d|, x(\delta_s + \varepsilon_d) > 0, y = 1\}.$$

Simplify the notation for products of indicator functions by writing, for example, $AH_+(z)$ for $A(z)H_+(z)$, and derive that

$$g_{a,b}(z) - g_{a,0}(z) - b'\Delta_2(z)$$
$$= \delta_d^2 + \varepsilon_s^2 + 2(\delta_s\delta_d + \varepsilon_s\varepsilon_d)[H_-(z) - H_+(z)]$$
(8)
$$+ 4(x\delta_d - \delta_s\delta_d - \varepsilon_s\varepsilon_d)[AH_-(z) - AH_+(z)]$$
$$+ 4(\delta_s - x + y\varepsilon_d)[A^0H_-(z) - A^0H_+(z)].$$

Define the remainder functions $s_{a,b}(z)$ by equalities

$$\|b\|s_{a,b}(z) = \delta_d^2 + \varepsilon_s^2 + 2(\delta_s\delta_d + \varepsilon_s\varepsilon_d)[H_-(z) - H_+(z)]$$
$$+ 4(x\delta_d - \delta_s\delta_d - \varepsilon_s\varepsilon_d)[AH_-(z) - AH_+(z)],$$

and observe that Property 1 holds for the class $\{s_{a,b}\}$. Thus, conditions (v) and (vi) of Theorem 3 are satisfied if the functions $l_{a,b}(z)$ are defined as the remaining part of expression (8), namely, $4(\delta_s - x + y\varepsilon_d)[A^0H_-(z) - A^0H_+(z)]$. Define $\tau_a = 1/4$ and $\tau_b = 1/2$ as Theorem 3 prescribes. It is only left to check condition (vii) of Theorem 3 by establishing

(9) $$\sup_{(a,b)\in\mathcal{N}_n} |\nu_n l_{a,b}| = o_p(n^{-1/2})$$

for all sequences of rectangles $\mathcal{N}_n$ of the order $O(n^{-1/4}) \times O(n^{-3/8})$ that are centered at the origin. Write out the Taylor approximations $U = \delta_s + \varepsilon_d + \delta_d\varepsilon_d - \varepsilon_s\varepsilon_d + o(\|(a,b)\|^2)$ and $D = \delta_s - \varepsilon_d - \delta_d\varepsilon_d - \varepsilon_s\varepsilon_d + o(\|(a,b)\|^2)$, and conclude that quantities $|D - (\delta_s - \varepsilon_d)|$ and $|U - (\delta_s + \varepsilon_d)|$ are of order $O(n^{-5/8})$ uniformly over $(a,b)$ in the neighborhoods $\mathcal{N}_n$. Use the oscillation properties of the empirical process established on page 765 in Shorack and Wellner [17] to conclude that

$$\sup_{(a,b)\in\mathcal{N}_n} |\nu_n AH_-(z) - \nu_n AH_+(z)| = o_p(n^{-1/4}).$$

Stochastic bound (9) follows directly.

Apply Theorem 3 and deduce that $(n^{1/4}a_n, n^{1/2}b_n) \rightsquigarrow (s^*, t^*)$, where
$$s^* = \arg\min_s[\psi_1(s) + s'Z_1]$$

and
$$t^* = \arg\min_t\left[\psi_2(t) + t'Z_2 + \sum_{i=1}^{3}\phi_i(s^*, t)\right].$$

A closed form expression for $(s^*, t^*)$ is given in Section 2.3 of Radchenko [13].

**8. Example: partial splines.** The following semiparametric example is discussed in Van de Geer ([19], Chapter 11), where a CLT is established for the parametric component using its characterization as a zero of the derivative of the criterion function. Below, the same result is derived by working directly with the definition of the estimator as a minimizer, using the approach introduced in Sections 2 and 3 for mixed-rates parametric problems.

Let $(Y_1, Z_1), \ldots, (Y_n, Z_n), \ldots$ be independent copies of $(Y, Z)$, where $Y$ is a real-valued response variable and $Z$ is a covariate. Suppose, for simplicity, that $Z$ takes values in $[0,1]^2$, write $Z = (U, V)$ and assume that the model

$$Y = g(U, V) + W$$

is satisfied with $E(W|Z) = 0$ and $g(U, V) = \theta U + \gamma(V)$. Here $\theta \in \mathbb{R}$ is an unknown parameter, and $\gamma$ is an unknown member of the functional class

$$\mathcal{S} = \left\{ \eta : [0,1] \to \mathbb{R}, \int_0^1 |\eta^{(m)}(v)|^2 \, dv < \infty \right\},$$

defined for a fixed positive integer $m$. Assume that the tails of the error distribution decrease exponentially fast: there exist positive constants $K$ and $\sigma_0^2$, such that, for all $z \in [0,1]^2$,

$$2K^2 E(e^{|W|/K} - 1 - |W|/K | Z = z) \leq \sigma_0^2.$$

Denote the distribution of $(U, V)$ by $Q$ and write $\|f\|_2$ for the $L_2(Q)$-norm of a function $f$. Define functions $e(v) = E(U|V = v)$ and $h(u, v) = u - e(v)$. Assume that $\|h\|_2 > 0$.

Fix a positive $\lambda_0$ and take $\lambda_n = \lambda_0 n^{-m/(2m+1)}$. Consider a class $\mathcal{F}$ of all regression functions $f$ of the form $f(u, v) = \alpha u + \eta(v)$ with $\alpha \in \mathbb{R}$ and $\eta \in \mathcal{S}$. Denote the roughness of such a function by $I^2(f) = I^2(\eta) = \int_0^1 |\eta^{(m)}(v)|^2 \, dv$. Define

$$g_n = \arg\min_{f \in \mathcal{F}} \left\{ \frac{1}{n} \sum_{i=1}^n [Y_i - f(U_i, V_i)]^2 + \lambda_n^2 I^2(f) \right\},$$

the penalized least squares estimator of function $g$ over the class $\mathcal{F}$. Assume that the regression of $U$ on $V$ is sufficiently smooth by requiring $I(e) < \infty$. Given a function $\tau$ from the class $\mathcal{S}$ and a real $\delta$, define $f_{\tau,\delta}(u, v) = [\theta + \delta]u + [\gamma(v) + \tau(v) - \delta e(v)]$ and note that function $f_{\tau,\delta}$ is a member of the class $\mathcal{F}$. Introduce criterion functions

$$G_n(\tau, \delta) = \frac{1}{n} \sum_{i=1}^n [Y_i - f_{\tau,\delta}(U_i, V_i)]^2 + \lambda_n^2 I^2(f_{\tau,\delta}).$$

Write $g_n(u, v)$ as $[\theta + \delta_n]u + [\gamma(v) + \tau_n(v) - \delta_n e(v)]$ and observe that the pair $(\tau_n, \delta_n)$ minimizes $G_n$ over the class $\{(\tau, \delta) : \tau \in \mathcal{S}, \delta \in \mathbb{R}\}$.

Methods from penalized least-squares estimation establish the common rate of convergence for the two components of the estimator $(\tau_n, \delta_n)$. Define $r = m/(2m+1)$. Stochastic bound $\|g_n - g\|_2 = O_p(n^{-r})$ is derived in Lemma 11.1 of Van de Geer[19]. Note that $\|g_n - g\|_2^2 = \delta_n^2 \|h\|_2^2 + \|\tau_n\|_2^2$, because conditional expectation $E(h(U,V)|V=v)$ is zero for each $v$ in $[0,1]$. Conclude that $\delta_n = O_p(n^{-r})$ and $\|\tau_n\|_2 = O_p(n^{-r})$.

Apply the approach of Section 3 to improve the convergence rate of $\delta_n$. Write $X_n$ for the standardized sum $n^{-1/2} \sum_{i=1}^n h(U_i, V_i) W_i$ and deduce that

$$
\begin{aligned}
(10) \quad & G_n(\tau, \delta) - G_n(\tau, 0) \\
& = \delta^2 Q_n h^2 - 2\delta[n^{-1/2} X_n - Q_n h \tau] + \lambda_n^2 [I^2(f_{\tau,\delta}) - I^2(f_{\tau,0})].
\end{aligned}
$$

Equality $Eh(U,V)\tau_n(V) = 0$ implies $Q_n h \tau_n = n^{-1/2} \nu_n h \tau_n$, and asymptotic equicontinuity of the empirical process indexed by functions $\{h\tau : \tau \in \mathcal{S}\}$ yields $Q_n h \tau_n = o_p(n^{-1/2})$. Use the definition of the roughness to derive $|I^2(f_{\tau_n,\delta}) - I^2(f_{\tau_n,0})| \le I(\delta e) I(2\gamma + 2\tau_n - \delta e)$. Note the stochastic bound $I(\tau_n) = O_p(1)$, implied by Van de Geer's Lemma 11.1, and conclude that $\lambda_n^2 [I^2(f_{\tau_n,\delta}) - I^2(f_{\tau_n,0})] = o_p(n^{-1/2}\delta)$. Expression (10) evaluated at $\tau = \tau_n$ and $\delta = \delta_n$ simplifies to

$$G_n(\tau_n, \delta_n) - G_n(\tau_n, 0) = \delta_n^2 Q_n h^2 - 2\delta_n n^{-1/2}[X_n + o_p(1)].$$

The law of large numbers yields $Q_n h^2 \to \|h\|_2$, and the limit is positive by assumption. Observe that $X_n = O_p(1)$ and apply Theorem 2, with $\delta^2 Q_n h^2$ playing the role of $M_n(\tau, \delta)$, to derive the correct $n^{-1/2}$ convergence rate of $\delta_n$.

Note that $\delta_n$ minimizes the criterion function $G_n(\tau_n, \delta) - G_n(\tau_n, 0)$ over $\delta$. Localize this function by writing $\delta = n^{-1/2} t$, and use the results of the previous paragraph to derive a quadratic approximation that holds uniformly on compacta,

$$(11) \quad G_n(\tau_n, n^{-1/2}t) - G_n(\tau_n, 0) = n^{-1}[t^2 Q_n h^2 - 2t X_n + o_p(1)].$$

Define $\sigma^2(z) = E(W^2|Z=z)$ and note that $X_n \rightsquigarrow X$, where $X \sim N(0, \|\sigma h\|_2^2)$. Minimization of the random quadratic function in (11) yields $n^{1/2}\delta_n = X_n/Q_n h^2 + o_p(1)$, and a CLT for $\delta_n$ follows directly. Note that because the criterion functions $G_n(\tau_n, \cdot)$ are convex, the formal derivation of the bound $\delta_n = O_p(n^{-1/2})$ could have been sidestepped.

## 9. Proofs.

9.1. *Proof of Theorem* 1. The next result is a version of the continuous mapping theorem.

LEMMA 2 (Modified continuous mapping). *Consider a metric space $(\mathcal{X}, d)$. Let random maps $X_n : A_n \to \mathcal{X}$ be defined on some sets $A_n \subset \Omega$ and consider a function $g : \mathcal{X} \to \mathbb{R}^d$ that is continuous at every point of a set $\mathcal{X}_0 \subset \mathcal{X}$. Suppose that $X : \Omega \to \mathcal{X}$ is a Borel measurable map for which there exists a Borel measurable set $A$, containing each of the sets $A_n$, such that $X \in \mathcal{X}_0$ on $A$. Suppose that $\mathbb{P}^*\{d(X_n, X) > \varepsilon\} \cap A_n \to 0$ for all $\varepsilon > 0$. Then $\mathbb{P}^*\{\|g(X_n) - g(X)\| > \delta\} \cap A_n \to 0$ for all $\delta > 0$.*

PROOF. Apply a standard device for proving continuous mapping theorems (see, e.g., the proof of Theorem 1.9.5 in van der Vaart and Wellner [21]). Fix a positive $\varepsilon$. Let $D_k$ be the set of all $x$ in $\mathcal{X}$ for which there exist $y$ and $z$ within the open ball of radius $1/k$ around $x$ with $\|g(y) - g(z)\| > \delta$. Note that $D_k$ is open and the sequence $D_k$ is decreasing. Also note that $\mathbb{P}\{X \in D_k\} \cap A \downarrow 0$, because every point in $\bigcap_{k=1}^\infty D_k$ is a point in $\mathcal{X}_0^c$. Observe that, for every fixed $k$,

$$\mathbb{P}^*\{\|g(X_n) - g(X)\| > \delta\} \cap A_n$$
$$\leq \mathbb{P}\{X \in D_k\} \cap A + \mathbb{P}^*\{d(X_n, X) > 1/k\} \cap A_n.$$

The first term on the right-hand side can be made arbitrarily small by choosing $k$ large enough. For a given choice of $k$, the second term tends to zero as $n$ goes to infinity. $\square$

Dudley [2] proved a representation theorem for the convergence in distribution in the sense of Hoffmann–Jørgensen. The following argument uses Dudley's result in the convenient form of Theorem 9.4 in Pollard [11], referred to as Representation Theorem.

PROOF OF THEOREM 1. Redefine function $H_n$ so that vector $(s_n, t_n)$ becomes the unique minimum. This can be done by leaving function $f_n$ unchanged and decreasing function $g_n$ by a $o_p^*(1)$ amount at exactly the point $(s_n, t_n)$. Note that the assumptions of the theorem remain valid after the change.

It is enough to show that $\mathbb{P}^* h(s_n, t_n) \to \mathbb{P}^* h(s^*, t^*)$ for all bounded, uniformly continuous, real functions $h$ on $\mathbb{R}^{d_1} \times \mathbb{R}^{d_2}$. Invoke the Representation Theorem for the convergence $(f_n, g_n) \rightsquigarrow (f, g)$, denote the corresponding perfect maps by $\phi_n$ and write $\widetilde{\omega}$ for the elements of the new probability space. Simplify the notation by replacing the composition $f_n(\phi_n(\widetilde{\omega}), s)$ with $\widetilde{f}_n(s)$, writing $\widetilde{s}_n$ for $s_n(\phi_n(\widetilde{\omega}))$, and so on, omitting the $\widetilde{\omega}$. Perfectness of $\phi_n$ implies $|\mathbb{P}^* h(s_n, t_n) - \mathbb{P} h(s^*, t^*)| \leq \widetilde{\mathbb{P}}^* |h(\widetilde{s}_n, \widetilde{t}_n) - h(\widetilde{s}^*, \widetilde{t}^*)|$, hence, it is enough to show that random vectors $(\widetilde{s}_n, \widetilde{t}_n)$ converge to $(\widetilde{s}^*, \widetilde{t}^*)$ in outer probability (see, e.g., Theorem 1.9.5 in van der Vaart and Wellner). Write $A_n^r$ for the subset of the new probability space that is defined by $\{\widetilde{\omega} : \|\widetilde{s}_n\| \vee \|\widetilde{t}_n\| \vee \|\widetilde{s}^*\| \vee \|\widetilde{t}^*\| \leq r\}$. Because quantities $\|\widetilde{s}_n\|, \|\widetilde{t}_n\|, \|\widetilde{s}^*\|$ and $\|\widetilde{t}^*\|$ are stochastically bounded with respect to $\widetilde{P}^*$, it is sufficient to check that, for each fixed $r$,

$$
\begin{aligned}
\widetilde{\mathbb{P}}^*\{\|\widetilde{s}_n - \widetilde{s}^*\| > \delta\} \cap A_n^r \to 0 \quad \text{and} \\
\widetilde{\mathbb{P}}^*\{\|\widetilde{t}_n - \widetilde{t}^*\| > \delta\} \cap A_n^r \to 0 \qquad \text{for all } \delta.
\end{aligned}
\tag{12}
$$

Fix a positive $r$ and restrict all the functions on $\mathbb{R}^{d_1}$ to the ball $\{\|x\| \leq r\}$. Denote by $\mathcal{X}_r$ the set of those (restricted) functions that are bounded and possess a unique minimum. Write $\Psi_r$ for the argmin map on $\mathcal{X}_r$. On the set $A_n^r$, the points $\widetilde{s}_n$ and $\widetilde{s}^*$ can be viewed as two values of the same map: $\widetilde{s}_n = \Psi_r[\alpha_n^{-1}\widetilde{H}_n(\cdot, \widetilde{t}_n)]$ and $\widetilde{s}^* = \Psi_r[\widetilde{f}]$. The two corresponding arguments are close when $n$ is large. Indeed, on the set $A_n^r$,

$$\sup_{\|s\| \leq r} |\alpha_n^{-1}\widetilde{H}_n(s, \widetilde{t}_n) - \widetilde{f}(s)|$$
$$\leq \sup_{\|s\| \leq r} |\widetilde{f}_n(s) - \widetilde{f}(s)| + \beta_n/\alpha_n \sup_{\|s\| \wedge \|t\| \leq r} |\widetilde{g}_n(s,t)|.$$

The right-hand side of the above inequality goes to zero in outer probability because of the bound from the Representation Theorem and the boundedness of $g(s,t)$. Let $A^r$ stand for the Borel measurable set $\{\widetilde{\omega} : \|\widetilde{s}^*\| \leq r\}$ and observe that on the set $A^r$ the map $\Psi_r$ is continuous at $\widetilde{f}$ [compare with condition (1) in Section 2]. Apply the modified continuous mapping lemma with $A_n^r$, $A^r$ and $\mathcal{X}_r$ playing the role of $A_n$, $A$ and $\mathcal{X}$, and deduce the first convergence in display (12).

Define $\Phi_r$ analogously to $\Psi_r$, but with respect to $\mathbb{R}^{d_2}$. Note that $\widetilde{t}_n = \Phi_r[\widetilde{g}_n(\widetilde{s}_n, \cdot)]$ and $\widetilde{t}^* = \Phi_r[\widetilde{g}(\widetilde{s}^*, \cdot)]$ on the set $A_n^r$. Also, on this set,

$$\sup_{\|t\| \leq r} |\widetilde{g}_n(\widetilde{s}_n, t) - \widetilde{g}(\widetilde{s}^*, t)|$$
$$\leq \sup_{\|s\| \wedge \|t\| \leq r} |\widetilde{g}_n(s,t) - \widetilde{g}(s,t)| + \sup_{\|t\| \leq r} |\widetilde{g}(\widetilde{s}_n, t) - \widetilde{g}(\widetilde{s}^*, t)|.$$

The first term on the right-hand side tends to zero in outer probability because of the bound from the Representation Theorem; the second term tends to zero in outer probability by the standard continuous mapping theorem. Deduce the second convergence in display (12) using an argument analogous to the one concluding the previous paragraph. □

9.2. *Proofs of the results in Section* 3. The following lemma simplifies the work with random polynomial functions.

LEMMA 3. *Let $\alpha$ be positive and let $\{\gamma_1, \ldots, \gamma_p, \eta_1, \ldots, \eta_p\}$ be a collection of nonnegative numbers satisfying $\gamma_i < \alpha$ for all $i \in \{1, \ldots, p\}$. Define* $\tau = \min_{i \leq p}(\frac{\eta_i}{\alpha - \gamma_i})$. *For each positive $\delta$ and each $O_p^*(1)$ sequence of random variables $L_n$, there exist a $O_p^*(1)$ sequence of random variables $M_n$, such that the following upper bound holds for all positive $u$:*

$$L_n\left(\sum_{i \leq p} n^{-\eta_i} u^{\gamma_i}\right) \leq \delta u^\alpha + M_n n^{-\alpha \tau}.$$

PROOF. It is enough to establish the bound for $\delta = 1$. Let $M_n(\omega)$ be the smallest real number satisfying the inequality $\sup_{u \geq 0}(L_n(\omega) \sum_{i \leq p} n^{-\eta_i} u^{\gamma_i} - u^\alpha) \leq M_n(\omega) n^{-\alpha \tau}$. Given a positive $\varepsilon$, select a large enough $L$ to ensure that $P^*\{L_n > L\} < \varepsilon$. Note that

$$P^*\{M_n > M\} \leq P^*\left\{\sup_{u \geq 0}\left(L \sum_{i \leq p} n^{-\eta_i} u^{\gamma_i} - u^\alpha\right) > M n^{-\alpha \tau}\right\} + \varepsilon.$$

To see that the first term on the right-hand side of the above inequality is zero for all $M$ large enough, combine the upper bound

$$\sup_{u \geq 0}\left(L \sum_{i \leq p} n^{-\eta_i} u^{\gamma_i} - u^\alpha\right) \leq \max_{i \leq k} \sup_{u \geq 0}(pL n^{-\eta_i} u^{\gamma_i} - u^\alpha)$$

with the inequalities

$$\sup_{u \geq 0}(pL n^{-\eta_i} u^{\gamma_i} - u^\alpha) = c_i n^{-\alpha \eta_i / (\alpha - \gamma_i)} \leq c_i n^{-\alpha \tau}, \qquad i = 1, \ldots, p.$$

Conclude that $M_n = O_p^*(1)$. □

Proof of Theorem 2 is omitted because it is similar to the following argument.

PROOF OF LEMMA 1. Deduce $\|a_n\|^\alpha + \|b_n\|^\beta = O_p^*(\sum_{i \leq p} n^{-\eta_i} \|(a_n, b_n)\|^{\gamma_i})$ from inequality $G_n(\alpha_n, b_n) - G_n(0, 0) \leq 0$. For each positive $\delta$, use Lemma 3 to establish

$$\|a_n\|^\alpha + \|b_n\|^\beta \leq \delta \|(a_n, b_n)\|^\alpha + O_p^*(n^{-\alpha \tau_a}).$$

Take a small enough $\delta$ and use inequality $\alpha \geq \beta$ to derive $\|a_n\|^\alpha + \|b_n\|^\beta = O_p^*(n^{-\alpha \tau_a})$. Conclude that $\|a_n\| = O_p^*(n^{-\tau_a})$ and $\|b_n\| = O_p^*(n^{-\alpha \tau_a / \beta})$. □

### 9.3. Proof of Theorem 3.

According to condition (iv),

$$G_n(a, b) = G(a, b) + n^{-1/2} a' \nu_n \Delta_1$$
$$+ n^{-1/2} b' \nu_n \Delta_2 + n^{-1/2} \|(a, b)\| \nu_n r_{a,b}.$$

Combine this representation with the stochastic bound $\|(a_n, b_n)\| = o_p(1)$ and deduce the approximation $G_n(a_n, b_n) = G(a_n, b_n) + O_p(n^{-1/2} \|(a_n, b_n)\|)$.

Apply Lemma 1 with function $G$ playing the role of $M_n$ and derive the stochastic bounds $\|a_n\| = O_p(n^{-\tau_a})$ and $\|b_n\| = O_p(n^{-\alpha\tau_a/\beta})$. It follows from condition (v) that

$$G_n(a,b) - G_n(a,0)$$
$$= G(a,b) - G(a,0) + n^{-1/2}b'\nu_n\Delta_2$$
$$+ n^{-1/2}\|b\|\nu_n s_{a,b} + n^{-1/2}\nu_n l_{a,b}.$$

Conditions (vi) and (vii) yield $\nu_n s_{a_n,b_n} = o_p(1)$ and $\nu_n l_{a_n,b_n} = o_p(n^{-\beta\tau_b+1/2})$. Thus,

(13)
$$G_n(a_n,b_n) - G_n(a_n,0)$$
$$= G(a_n,b_n) - G(a_n,0) + O_p(n^{-1/2}\|b_n\|) + o_p(n^{-\beta\tau_b}).$$

Observe that $\sum_{i=1}^{\alpha}\|a_n\|^{\alpha-i}\|b_n\|^i = O_p(n^{-1/2}\|b_n\|)$. Consequently,

$$G(a_n,b_n) - G(a_n,0)$$
$$= \psi_2(b_n)[1 + o_p(1)] + O_p\left(\sum_{i=1}^{p} n^{-\tau_a\gamma_i}\|b_n\|^{\eta_i}\right) + o_p(n^{-1/2}\|b_n\|).$$

Combine this approximation with approximation (13) and let the term $\psi_2(b_n)[1 + o_p(1)]$ play the role of $M_n(a_n,b_n)$ in Theorem 2. Conclude that $\|b_n\| = O_p(n^{-\tau_b})$.

Introduce new variables $s$ and $t$ by $s = n^{\tau_a}a$ and $t = n^{\tau_b}b$. Observe that

$$\sum_{i=1}^{p}\phi_i(n^{-\tau_a}s, n^{-\tau_b}t) = n^{-\beta\tau_b}\sum_{i=1}^{p} n^{-(\beta-\eta_i)[\lambda_i-\tau_b]}\phi_i(s,t)$$

and

$$\sum_{i=1}^{\alpha}\|n^{-\tau_a}s\|^{\alpha-i}\|n^{-\tau_b}t\|^i \leq n^{-\tau_a(\alpha-1)-\tau_b}\sum_{i=1}^{\alpha}\|s\|^{\alpha-i}\|t\|^i$$
$$= n^{-\beta\tau_b}n^{-(\beta-1)[\lambda_0-\tau_b]}\sum_{i=1}^{\alpha}\|s\|^{\alpha-i}\|t\|^i.$$

Combine the last two displays with the approximations in conditions (iv) and (vii) to deduce

(14)
$$G(n^{-\tau_a}s, n^{-\tau_b}t) = n^{-\alpha\tau_a}[\psi_1(s) + q_n(s)]$$
$$+ n^{-\beta\tau_b}[\psi_2(t) + \phi_i(s,t)1\{\lambda_i = \tau_b\} + w_n(s,t)],$$

where $\sup_{K_1}|q_n(s)| = o(1)$ for every compact set $K_1$ in $\mathbb{R}^{d_1}$, and $\sup_{K_2}|w_n(s,t)| = o(1)$ for every compact set $K_2$ in $\mathbb{R}^{d_1} \times \mathbb{R}^{d_2}$. Conditions (iv) through

(vii) yield

$$G_n(n^{-\tau_a}s, n^{-\tau_b}t) = G(n^{-\tau_a}s, n^{-\tau_b}t) + n^{-\alpha\tau_a}[s'\nu_n\Delta_1 + q'_n(s)] \quad (15)$$
$$+ n^{-\beta\tau_b}[n^{-(\beta-1)[\lambda_0-\tau_b]}t'\nu_n\Delta_2 + w'_n(s,t)],$$

where $\sup_{K_1} |q'_n(s)| = o_p(1)$ for every compact set $K_1$ in $\mathbb{R}^{d_1}$, and $\sup_{K_2} |w'_n(s,t)| = o_p(1)$ for every compact set $K_2$ in $\mathbb{R}^{d_1} \times \mathbb{R}^{d_2}$.

Denote $G_n(n^{-\tau_a}s, n^{-\tau_b}t)$ by $H_n(s,t)$. Combine approximations (14) and (15) and conclude that, uniformly on compacta in $\mathbb{R}^{d_1} \times \mathbb{R}^{d_2}$,

$$H_n(s,t) = n^{-\alpha\tau_a}[\psi_1(s) + s'\nu_n\Delta_1 + o_p(1)]$$
$$+ n^{-\beta\tau_b}\left[\psi_2(t) + 1\{\lambda_0 = \tau_b\}t'\nu_n\Delta_2 \right.$$
$$\left. + \sum_{i=1}^{p} 1\{\lambda_i = \tau_b\}\phi_i(s,t) + o_p(1)\right].$$

Note that $\alpha\tau_a \leq \beta\tau_b$. Indeed, this inequality is valid in the case $\lambda_0 = \tau_b$; in the case $\lambda_0 \neq \tau_b$, it follows from approximation (14) and the fact that function $G$ assumes only nonnegative values near the origin. If $\alpha\tau_a < \beta\tau_b$, apply Theorem 1 to complete the proof. If $\alpha\tau_a = \beta\tau_b$, the standard arg min theorem will suffice.

**Acknowledgments.** I thank David Pollard, the referees and the Associate Editor for their valuable suggestions that led to an improved presentation of the paper.

UNIVERSITY OF SOUTHERN CALIFORNIA
3670 TROUSDALE PKWY, BRIDGE HALL 401-M
LOS ANGELES, CALIFORNIA 90089-0809
USA
E-MAIL: radchenk@usc.edu
URL: http://www-rcf.usc.edu/~radchenk/